\documentclass{amsart}

\newtheorem{theorem}{Theorem}
\newtheorem{lemma}[theorem]{Lemma}
\newtheorem{proposition}[theorem]{Proposition}
\newtheorem{corollary}[theorem]{Corollary}

\theoremstyle{definition}

\newtheorem{remark}[theorem]{Remark}

\usepackage{amscd,amssymb}

\begin{document}

\title[Invariants of unipotent transformations]
{Invariants of unipotent transformations\\
acting on noetherian relatively free algebras}

\author[Vesselin Drensky]
{Vesselin Drensky}
\address{Institute of Mathematics and Informatics,
Bulgarian Academy of Sciences,
1113 Sofia, Bulgaria}
\email{drensky@math.bas.bg}

\thanks
{Partially supported by Grant MM-1106/2001 of the Bulgarian National Science Fund.}


\begin{abstract} 
The classical theorem of Weitzenb\"ock states that the algebra 
of invariants $K[X]^g$
of a single unipotent transformation $g\in GL_m(K)$ 
acting on the polynomial algebra $K[X]=K[x_1,\ldots,x_m]$
over a field $K$ of characteristic 0
is finitely generated. This algebra coincides with the algebra of constants
$K[X]^{\delta}$ of a linear locally nilpotent derivation $\delta$ of $K[X]$.
Recently the author and C.~K. Gupta have started the study
of the algebra of invariants $F_m({\mathfrak V})^g$
where $F_m({\mathfrak V})$ is the relatively free algebra of rank $m$ in a variety 
$\mathfrak V$ of associative algebras. They have shown that
$F_m({\mathfrak V})^g$ is not finitely generated if $\mathfrak V$ contains
the algebra $UT_2(K)$ of $2\times 2$ upper triangular matrices.
The main result of the present paper is that
the algebra $F_m({\mathfrak V})^g$ is finitely generated if 
and only if the variety ${\mathfrak V}$ does not contain the algebra $UT_2(K)$.
As a by-product of the proof we have established also the finite generation
of the algebra of invariants $T_{nm}^g$ where $T_{nm}$ is the mixed trace algebra
generated by $m$ generic $n\times n$ matrices and the traces of their products.
\end{abstract}

\maketitle

\section*{Introduction}

Let $K$ be any field of characteristic 0 and let $X=\{x_1,\ldots,x_m\}$,
where $m>1$. Let $g\in GL_m=GL_m(K)$ be a unipotent
linear operator acting on the vector space 
$KX=Kx_1\oplus\cdots\oplus Kx_m$. 
By the classical theorem of Weitzenb\"ock \cite{W}, the algebra of invariants
\[
K[X]^g=\{f\in K[X]\mid f(g(x_1),\ldots,g(x_m))=f(x_1,\ldots,x_m)\}
\]
is finitely generated. A proof in modern language was given by Seshadri \cite{S}.
An elementary proof based on the ideas of \cite{S}
was presented by Tyc \cite{T}.
Since $g-1$ is a nilpotent linear operator of $KX$, we may consider
the linear locally nilpotent derivation 
\[
\delta=\log g=\sum_{i\geq 1}(-1)^{i-1}\frac{(g-1)^i}{i}
\]
called a Weitzenb\"ock derivation. 
(The $K$-linear operator $\delta$ acting on an algebra $A$ is called a 
derivation if $\delta(uv)=\delta(u)v+u\delta(v)$ for all $u,v\in A$.)
The algebra of invariants
${\mathbb C}[X]^g$ coincides with the algebra of constants 
${\mathbb C}[X]^{\delta}\ (=\ker(\delta)$~). See the book by Nowicki \cite{N}
for a background on the properties of the algebras of constants
of Weitzenb\"ock derivations.

Looking for noncommutative generalizations of invariant theory,
see e.~g. the survey by Formanek \cite{F1}, let
$K\langle X\rangle=K\langle x_1,\ldots,x_m\rangle$ be the free 
unitary associative algebra freely generated by $X$. 
The action of $GL_m$
is extended diagonally on $K\langle X\rangle$ by the rule
\[
h(x_{j_1}\cdots x_{j_n})=h(x_{j_1})\cdots h(x_{j_n}),\
h\in GL_m,\ x_{j_1},\ldots,x_{j_n}\in X.
\]
For any PI-algebra $R$, let $T(R)\subset K\langle X\rangle$ 
be the T-ideal of all polynomial identities in $m$ variables
satisfied by $R$. The class $\mathfrak V=\text{var}(R)$ of all algebras
satisfying the identities of $R$ is called the variety of algebras generated by $R$
(or determined by the polynomial identities of $R$). The factor algebra
$F_m({\mathfrak V})=K\langle X\rangle/T(R)$ is called the relatively free
algebra of rank $m$ in $\mathfrak V$.
We shall use the same symbols $x_j$ and $X$ for the generators
of $F_m({\mathfrak V})$. Since $T(R)$ is $GL_m$-invariant, the action of $GL_m$
on $K\langle X\rangle$ is inherited by $F_m({\mathfrak V})$ and one can consider
the algebra of invariants $F_m({\mathfrak V})^G$ for any linear group $G$.
As in the commutative case, if $g\in GL_m$ is unipotent, then 
$F_m({\mathfrak V})^g$ coincides with the algebra 
$F_m({\mathfrak V})^{\delta}$ of the constants of the derivation $\delta=\log g$.

Till the end of the paper we fix the integer $m>1$, 
the variety $\mathfrak V$, the unipotent
linear operator $g\in GL_m$ and the derivation $\delta=\log g$.

The author and C.~K. Gupta \cite{DG} have started the study
of the algebra of invariants $F_m({\mathfrak V})^g$.
They have shown that if $\mathfrak V$ contains the algebra $UT_2(K)$ 
of $2\times 2$ upper triangular matrices and $g$ is different from the identity
of $GL_m$, then
$F_m({\mathfrak V})^g$ is not finitely generated for any $m>1$.
They have also established that, if $UT_2(K)$ does not belong to $\mathfrak V$,
then, for $m=2$, the algebra $F_2({\mathfrak V})^g$ is finitely generated. 

In the present paper we close the problem for which varieties $\mathfrak V$
and which $m$ the algebra $F_m({\mathfrak V})^g$ is finitely generated. 
Our main result is that this holds, and for all $m>1$, 
if and only if the variety ${\mathfrak V}$ does not contain the algebra $UT_2(K)$.

It is natural to expect such a result by two reasons. First, it follows from the proof
of Tyc \cite{T}, see also the earlier paper by Onoda \cite{O}, that 
the algebra $K[X]^g$ is isomorphic to the algebra of invariants 
of certain $SL_2$-action on the polynomial algebra in $m+2$ variables.
One can prove a similar fact for $F_m({\mathfrak V})^g$
and $(K[y_1,y_2]\otimes_KF_m({\mathfrak V}))^{SL_2}$. Second, 
the results of Vonessen \cite{V}, Domokos and the author \cite{DD} 
give that $F_m({\mathfrak V})^G$ is finitely
generated for all reductive $G$ if and only if the finitely generated algebras
in $\mathfrak V$ are one-side noetherian. For unitary algebras this means that
$\mathfrak V$ does not contain $UT_2(K)$ or, equaivalently,
$\mathfrak V$ satisfies the Engel identity $[x_2,x_1,\ldots,x_1]=0$.
In our proof we use the so called proper polynomial identities
introduced by Specht \cite{Sp}, the fact that the Engel identity
implies that the vector space of proper polynomials in $F_m({\mathfrak V})$
is finite dimensional and hence $F_m({\mathfrak V})$ has a series of ideals
such that the factors are finitely generated $K[X]$-modules.
As a by-product of the proof we have established also the finite generation
of the algebra of invariants $T_{nm}^g$, where $T_{nm}$ is the mixed trace algebra
generated by $m$ generic $n\times n$ matrices $x_1,\ldots,x_m$
and and the traces of their products $\text{\rm tr}(x_{i_1}\cdots x_{i_k})$,
$k\geq 1$.

\section{Preliminaries}

We fix two finite dimensional vector spaces $U$ and $V$,
$\dim U=p$, $\dim V=q$, and representations of the infinite cyclic group
$G=\langle g\rangle$:
\[
\rho_U:G\to GL(U)=GL_p,\quad \rho_V:G\to GL(V)=GL_q,
\]
where $\rho_U(g)$ and $\rho_V(g)$ are unipotent linear operators. Fixing bases
$Y=\{y_1,\ldots,y_p\}$ and $Z=\{z_1,\ldots,z_q\}$ of $U$ and $V$, respectively,
we consider the free left $K[Y]$-module $M(Y,Z)$ with basis $Z$. Then $g$ acts 
diagonally on $M(Y,Z)$ by the rule
\[
g:\sum_{j=1}^qf_j(y_1,\ldots,y_p)z_j\to 
\sum_{j=1}^qf_j(g(y_1),\ldots,g(y_p))g(z_j),\quad f_j\in K[Y],
\]
where, by definition, $g(y_i)=\rho_U(g)(y_i)$ and $g(z_j)=\rho_V(g)(z_j)$. Let
$M(Y,Z)^g$ be the set of fixed points in $M(Y,Z)$ under the action of $g$.
Since $\rho_U(g)$ and $\rho_V(g)$ are unipotent operators, the operators
$\delta_U=\log\rho_U(g)$ and $\delta_V=\log\rho_V(g)$ are well defined. Denote
by $\delta$ the induced derivation of $K[Y]$. We extend $\delta$ to 
a derivation of $M(Y,Z)$, denoted also by $\delta$, i.~e. $\delta$ is the linear 
operator of $M(Y,Z)$ defined by
\[
\delta:\sum_{j=1}^qf_j(Y)z_j\to \sum_{j=1}^q\delta(f_j(Y))z_j
+\sum_{j=1}^qf_j(Y)\delta(z_j).
\]
It is easy to see that $\delta=\log g$ on $M(Y,Z)$ and $M(Y,Z)^g$ coincides with
the kernel of $\delta$, i.~e. the set of constants $M(Y,Z)^{\delta}$.
Changing the bases of $U$ and $V$, we may assume that $\delta_U$ and $\delta_V$ 
have the form
\[
\delta_U=\begin{pmatrix}
J_{p_1}&0&\cdots&0&0\\
0&J_{p_2}&\cdots&0&0\\
\vdots&\vdots&\cdots&\vdots&\vdots\\
0&0&\cdots&J_{p_{k-1}}&0\\
0&0&\cdots&0&J_{p_k}\\
\end{pmatrix},\quad\ 
\delta_V=\begin{pmatrix}
J_{q_1}&0&\cdots&0&0\\
0&J_{q_2}&\cdots&0&0\\
\vdots&\vdots&\cdots&\vdots&\vdots\\
0&0&\cdots&J_{q_{l-1}}&0\\
0&0&\cdots&0&J_{q_l}\\
\end{pmatrix},
\]
where $J_r$ is the $(r+1)\times(r+1)$ Jordan cell 
\begin{equation}\label{Jordan cell}
J_r=\begin{pmatrix}
0&1&0&\cdots&0&0\\
0&0&1&\cdots&0&0\\
\vdots&\vdots&\vdots&\cdots&\vdots&\vdots\\
0&0&0&\cdots&1&0\\
0&0&0&\cdots&0&1\\
0&0&0&\cdots&0&0\\
\end{pmatrix}
\end{equation}
with zero diagonal.

We denote by $W_r$ the irreducible $(r+1)$-dimensional $SL_2$-module.
It is isomorphic to the $SL_2$-module of the forms of degree $r$ in two variables
$x,y$. This is the unique structure of an $SL_2$-module on the $(r+1)$-dimensional vector 
space which agrees with the action of $\delta$ (and hence of $g$) as the Jordan cell
(\ref{Jordan cell}): We can think of $\delta$ as the derivation of $K[x,y]$ defined
by $\delta(x)=0$, $\delta(y)=x$. We fix the
``canonical'' basis of $W_r$
\begin{equation}\label{module of forms in two variables}
u^{(0)}=x^r,u^{(1)}=\frac{x^{r-1}y}{1!},u^{(2)}=\frac{x^{r-2}y^2}{2!},\ldots,
u^{(r-1)}=\frac{xy^{r-1}}{(r-1)!},u^{(r)}=\frac{y^r}{r!}.
\end{equation}
We give $U$ and $V$ the structure of $SL_2$-modules 
\begin{equation}\label{U and V are SL2-modules}
U=W_{p_1}\oplus\cdots\oplus W_{p_k},\quad V=W_{q_1}\oplus\cdots\oplus W_{q_l},
\end{equation}
and extend it on $K[Y]$ and $M(Y,Z)$ via the diagonal action of $SL_2$.
Again, this agrees with the action of $g$ and $\delta$.
Then $K[U]$ and $M(Y,Z)$ are direct sums of irreducible $SL_2$-modules
$U_{ri}\subset K[Y]$ and $W_{rj}\subset M(Y,Z)$ isomorphic to $W_r$,
$i,j=1,2,\ldots$, $r=0,1,2,\ldots$, with canonical bases
$\{u_{ri}^{(0)},u_{ri}^{(1)},\ldots,u_{ri}^{(r)}\}$ and 
$\{w_{rj}^{(0)},w_{rj}^{(1)},\ldots,w_{rj}^{(r)}\}$, respectively.

\begin{lemma}\label{lemma 1}
The elements $u\in K[Y]$ and $w\in M(Y,Z)$ belong to $K[Y]^{\delta}$ and
$M(Y,Z)^{\delta}$, respectively, if and only if they have the form
\begin{equation}\label{canonical form of constants}
u=\sum_{r,i}\alpha_{ri}u_{ri}^{(0)},\quad
w=\sum_{r,j}\beta_{rj}w_{rj}^{(0)},\quad
\alpha_{ri},\beta_{rj}\in K.
\end{equation}
\end{lemma}

\begin{proof}
Almkvist, Dicks and Formanek \cite{ADF} translated in the language of $g$-invariants
results of De Concini, Eisenbud and Procesi \cite{DEP} and proved that, in our
notation, $g(u)=u$ and $g(w)=w$ if and only if $u$ and $w$ have the form
(\ref{canonical form of constants}). Since $g(u)=u$ if and only if $\delta(u)=0$,
and similarly for $w$, we obtain that (\ref{canonical form of constants})
holds if and only if $u$ and $w$ are $\delta$-constants. (The same fact is
contained in the paper by Tyc \cite{T} but in the language of representations
of the Lie algebra $sl_2(K)$.)
\end{proof}

In each component $W_r$ of $U$ in (\ref{U and V are SL2-modules}), using the basis
(\ref{module of forms in two variables}), we define a linear operator $d$ by
\[
d(u^{(k)})=(k+1)(r-k)u^{(k+1)},\quad k=0,1,2,\ldots,r,
\]
i.~e., up to multiplicative constants, $d$ acts by
$u^{(0)}\to u^{(1)}\to u^{(2)}\to\cdots\to u^{(r)}\to 0$. We extend $d$ to a derivation
of $K[Y]$. As in the case of $\delta$, again we can think of $d$ as the derivation
of $K[x,y]$ defined by $d(x)=y$, $d(y)=0$.

\begin{lemma}\label{lemma 2}
{\rm (i)} The derivation $d$ acts on each irreducible component $U_{ri}$ of $K[Y]$ by
\[
d(u_{ri}^{(k)})=(k+1)(r-k)u_{ri}^{(k+1)},\quad k=0,1,\ldots,r.
\]

{\rm (ii)} If $f=f(Y)\in K[Y]$, then $\delta^{s+1}(f)=0$ if and only if $f$
belongs to the vector space
\begin{equation}\label{constants of level s}
K[Y]_s=\sum_{t=0}^sd^t(K[Y]^{\delta}).
\end{equation}
\end{lemma}

\begin{proof} Part (i) follows from the fact that the $SL_2$-action
on $U$ is the only action which agrees with the action of $\delta$ as well as 
with the action of $d$ (as the derivations of $K[x,y]$ defined by
$\delta(x)=0$, $\delta(y)=x$ and $d(x)=y$, $d(y)=0$, respectively),
and the extension of this $SL_2$-action to $K[U]$ also agrees with
the action of $\delta$ and $d$ on $K[U]$. 

(ii) Since the irreducible $SL_2$-submodules of $K[Y]$ are $\delta$- and $d$-invariant,
it is sufficient to prove the statement only for $f\in W_r\subset K[Y]$.
Considering the basis (\ref{module of forms in two variables}) of $W_r$, we have that
$\delta^{s+1}(f)=0$ if and only if 
\[
f=\alpha_0u^{(0)}+\alpha_1u^{(1)}+\cdots+\alpha_su^{(s)},\quad \alpha_k\in K.
\]
Since $W_r^{\delta}=Ku^{(0)}$ and $d^t(u^{(0)})\in Ku^{(t)}$, we obtain that
$W_r\cap K[Y]_s$ is spanned by $u^{(0)},u^{(1)},\ldots,u^{(s)}$ and coincides with
the kernel of $\delta^{s+1}$ in $W_r$.
\end{proof}

In principle, the proof of the following proposition can be obtained following
the main steps of the proof of Tyc \cite{T} of the Weitzenb\"ock theorem.
The proof of the three main lemmas in \cite{T} uses only the fact that the ideals
of the algebra $K[Y]$ are finitely generated $K[Y]$-modules. Instead, we shall
give a direct proof, using the idea of the proof of Lemma 3 in \cite{T}.

\begin{proposition}\label{proposition 3}
The set of constants $M(Y,Z)^{\delta}$ is a finitely generated 
$K[Y]^{\delta}$-module.
\end{proposition}

\begin{proof}
Let $N_i$ be the $K[Y]$-submodule of $M(Y,Z)$ generated by the basis elements $z_j$ 
of $V=Kz_1\oplus\cdots\oplus Kz_q$ corresponding to the $i$-th Jordan cell $J_{q_i}$.
Since $M(Y,Z)=N_1\oplus\cdots\oplus N_l$ and 
$M(Y,Z)^{\delta}=N_1^{\delta}\oplus\cdots\oplus N_l^{\delta}$, it is sufficient to
show that each $N_i^{\delta}$ is a finitely generated $K[Y]^{\delta}$-module.
Hence, without loss of generality we may assume that $q=r+1$ and
$\delta(z_0)=0$, $\delta(z_j)=z_{j-1}$, $j=1,2,\ldots,r$. Let
\begin{equation}\label{expression of constants}
f=f_0(Y)z_0+f_1(Y)z_1+\cdots+f_r(Y)z_r\in M(Y,Z)^{\delta},\quad f_j(Y)\in K[Y].
\end{equation}
Then
\[
\delta(f)=(\delta(f_0)+f_1)z_0+(\delta(f_1)+f_2)z_1+\cdots
+(\delta(f_{r-1})+f_r)z_{r-1}+\delta(f_r)z_r
\]
and this implies that 
\[
\delta(f_j)=-f_{j+1},\quad j=0,1,\ldots,r-1,
\]
\[
\delta(f_r)=\delta^2(f_{r-1})=\cdots=\delta^r(f_1)=\delta^{r+1}(f_0)=0.
\]
Hence, fixing any element $f_0(Y)$ from $K[Y]_r$, we determine all the
coefficients $f_1,\ldots,f_r$ from (\ref{expression of constants}).
By Lemma \ref{lemma 2} it is sufficient to show that 
the $K[Y]^{\delta}$-module generated by $d^t(K[Y]^{\delta})$
is finitely generated. By the theorem of Weitzenb\"ock, $K[Y]^{\delta}$ is
a finitely generated algebra. Let $\{h_1,\ldots,h_n\}$ be a generating set of
$K[Y]^{\delta}$. Then $d^t(K[Y]^{\delta})$ is spanned by the elements
$d^t(h_1^{a_1}\cdots h_n^{a_n})$. Since $d$ is a derivation, 
$d^t(K[Y]^{\delta})$ is spanned by elements of the form
\[
h_1^{c_1}\cdots h_n^{c_n}\left(\prod d^{t_{i_1}}(h_1)\right)
\cdots \left(\prod d^{t_{i_n}}(h_n)\right),\quad
\sum t_{i_1}+\cdots+\sum t_{i_n}=t.
\]
There is only a finite number of possibilities for $t_{i_1},\ldots,t_{i_n}$,
and we obtain that $d^t(K[Y]^{\delta})$ generates a finitely generated
$K[Y]^{\delta}$-module.
\end{proof}

\begin{corollary}\label{corollary of proposition 3}
Let, in the notation of this section, $U$ and $V$ be polynomial $GL_m$-modules,
let $g\in GL_m$ be a unipotent matrix and let $M(Y,Z)$ be equipped with
the diagonal action of $GL_m$. Then, for every $GL_m$-submodule $M_0$
of $M(Y,Z)$, the natural homomorphism $M(Y,Z)\to M(Y,Z)/M_0$
induces an epimorphism $M(Y,Z)^g\to (M(Y,Z)/M_0)^g$,
i.~e. we can lift the $g$-invariants of $M(Y,Z)/M_0$ to $g$-invariants
of $M(Y,Z)$.
\end{corollary}

\begin{proof}
The lifting of the constants was established in \cite{DG} in the case of
relatively free algebras and the same proof works in our situation. 
Since $U$ and $V$ are polynomial $GL_m$-modules, the module $M(Y,Z)$ is
completely reducible. Hence $M(Y,Z)=M_0\oplus M'$ for some $GL_m$-submodule
$M'$ of $M(Y,Z)$ and $M(Y,Z)/M_0\cong M'$. If $w+M_0=\bar w\in (M(Y,Z)/M_0)^g$,
then $w=w_0+w'$, $w_0\in M_0$, $w'\in M'$, and
$g(w)=g(w_0)+g(w')$. Since $g(\bar w)=\bar w$, we obtain that $g(w')=w'$
and the $g$-invariant $\bar w$ is lifted to the $g$-invariant $w'$.
\end{proof}

\begin{remark}
The proof of Proposition \ref{proposition 3} gives also an algorithm to find the
generators of $M(Y,Z)^{\delta}$ in terms of the generators of $K[Y]^{\delta}$.
The latter problem is solved by van den Essen \cite{E} and his algorithm uses
Gr\"obner bases techniques.
\end{remark}

\section{The Main Results}

The following theorem is the main result of our paper. For $m=2$ it was established
in \cite{DG} using the description of the $g$-invariants of $K\langle x,y\rangle$.

\begin{theorem}\label{main theorem}
For any variety $\mathfrak V$ of associative algebras which does not contain the
algebra $UT_2(K)$ of $2\times 2$ upper triangular matrices, 
the algebra of invariants $F_m({\mathfrak V})^g$  of any unipotent
$g\in GL_m$ is finitely generated.
\end{theorem}

\begin{proof}
We shall work with the linear locally nilpotent derivation $\delta=\log g$
instead with $g$.

It is well known that any variety $\mathfrak V$ which does not contain
$UT_2(K)$ satisfies some Engel identity $[x_2,x_1,\ldots,x_1]=0$.
By the theorem of Zelmanov \cite{Z} any Lie algebra over 
a field of characteristic zero satisfying the Engel identity is nilpotent.
Hence we may assume that $\mathfrak V$ satisfies 
the polynomial identity of Lie nilpotency $[x_1,\ldots,x_{c+1}]=0$.
(Actually, this follows from much easier and much earlier results on PI-algebras.)

Let us consider the vector space $B_m({\mathfrak V})$
of so called proper polynomials in $F_m({\mathfrak V})$. It is spanned by all
products $[x_{i_1},\ldots,x_{i_k}]\cdots[x_{j_1},\ldots,x_{j_l}]$
of commutators of length $\geq 2$. 
One of the main results of the paper by the author \cite{D} states that
if $\{f_1,f_2,\ldots\}$ is a basis of $B_m({\mathfrak V})$, then
$F_m({\mathfrak V})$ has a basis 
\[
\{x_1^{p_1}\cdots x_m^{p_m}f_i\mid p_j\geq 0,
i=1,2,\ldots\}.
\]
Let $B_m^{(k)}({\mathfrak V})$
be the homogeneous component of degree $k$ of $B_m({\mathfrak V})$.
It follows from the proof of Theorem 5.5 in \cite{D},
that for any Lie nilpotent variety $\mathfrak V$, 
and for a fixed positive integer $m$,
the vector space $B_m({\mathfrak V})$ 
is finite dimensional. Hence $B_m^{(n)}({\mathfrak V})=0$ for $n$ sufficiently large, 
e.~g. for $n>n_0$. Let $I_k$ be the ideal of $F_m({\mathfrak V})$ generated by
$B_m^{(k+1)}({\mathfrak V})+B_m^{(k+2)}({\mathfrak V})+\cdots+B_m^{(n_0)}({\mathfrak V})$.
Since $wx_i=x_iw+[w,x_i]$, $w\in F_m({\mathfrak V})$, applying Lemma 2.4 \cite{D},
we obtain that $I_k/I_{k+1}$ is a free left $K[X]$-module with any basis of
the vector space $B_m^{(k)}({\mathfrak V})$ as a set of free generators.
Since $\delta$ is a nilpotent linear operator of $U=KX=Kx_1\oplus \cdots\oplus Kx_m$, 
it acts also as a nilpotent linear operator of $V_k=B_m^{(k)}({\mathfrak V})$.
Proposition \ref{proposition 3} gives that
$(I_k/I_{k+1})^{\delta}$ is a finitely generated $K[X]^{\delta}$-module.
Clearly, $B_m^{(0)}({\mathfrak V})=K$, $B_m^{(1)}({\mathfrak V})=0$,
$B_m^{(2)}({\mathfrak V})$ is spanned by the commutators $[x_{i_1},x_{i_2}]$, etc.
Hence $I_0/I_1\cong K[X]$ and by the theorem of Weitzenb\"ock $(I_0/I_1)^{\delta}$
is a finitely generated algebra. We fix a system of generators 
$\bar f_1,\ldots,\bar f_a$ of the algebra $(I_0/I_1)^{\delta}$ and finite sets of
generators $\{\bar f_{k1},\ldots,\bar f_{kb_k}\}$ of the $K[X]^{\delta}$-modules
$(I_k/I_{k+1})^{\delta}$, $k=2,3,\ldots,n_0$. 
The vector space $U$ is a $GL_m$-module and its $GL_m$-action
makes $V_k$ a polynomial $GL_m$-module. We apply Corollary 
\ref{corollary of proposition 3} and lift all $\bar f_i$ and $\bar f_{kj}$
to some $\delta$-constants $f_i,f_{kj}\in F_m({\mathfrak V})^{\delta}$.
The algebra $S$ generated by $f_1,\ldots,f_a$ maps onto $(I_0/I_1)^{\delta}$
and hence $(I_k/I_{k+1})^{\delta}$ is a left $S$-module generated by
$\bar f_{k1},\ldots,\bar f_{kb_k}$. The condition $I_{n_0+1}=0$ together with
Corollary \ref{corollary of proposition 3} gives that the $f_i$ and $f_{kj}$
generate $F_m({\mathfrak V})^{\delta}$.
\end{proof}

Together with the results of \cite{DG} Theorem \ref{main theorem} gives immediately:

\begin{corollary}
For $m\geq 2$ and for any fixed unipotent operator $g\in GL_m$,
$g\not=1$, the algebra of $g$-invariants $F_m({\mathfrak V})^g$
is finitely generated if and only if $\mathfrak V$ does not contain
the algebra $UT_2(K)$.
\end{corollary}

We refer to the books \cite{F2} and \cite{DF} for a background on the theory of matrix
invariants. We fix an integer $n>1$ and consider the generic $n\times n$ matrices
$x_1,\ldots,x_m$. Let $C_{nm}$ be the 
pure trace algebra, i.~e. the algebra generated by the traces of products 
$\text{tr}(x_{i_1}\cdots x_{i_k})$, $k=1,2,\ldots$, and let $T_{nm}$ be
the mixed trace algebra generated by $x_1,\ldots,x_m$ and $C_{nm}$.
It is well known that $C_{nm}$ is finitely generated. (The Nagata-Higman theorem
states that the nil polynomial identity $x^n=0$ implies the identity of nilpotency
$x_1\cdots x_d=0$. If $d(n)$ is the minimal $d$ with this property, one may take
as generators $\text{tr}(x_{i_1}\cdots x_{i_k})$ with $k\leq d(n)$.)
Also, $T_{nm}$ is a finitely generated $C_{nm}$-module.

\begin{theorem}
For any unipotent operator $g\in GL_m$, the algebra $T_{nm}^g$
is finitely generated.
\end{theorem}

\begin{proof}
Consider the vector space $U$ of all formal traces
$y_i=\text{tr}(x_{i_1}\cdots x_{i_k})$, $i_j=1,\ldots,m$, $1\leq k\leq d(n)$.
Let $Y$ be the set of all $y_i$.
It has a natural structure of a $GL_m$-module and hence $\delta=\log g$
acts as a nilpotent linear operator on $U$. Also, consider a finite system of
generators $f_1,\ldots,f_a$ of the $C_{nm}$-module $T_{nm}$. We may assume that the
$f_j$ do not depend on the traces and fix some elements $h_j\in K\langle X\rangle$
such that $h_j\to f_j$ under the natural homomorphism 
$K\langle X\rangle\to T_{nm}$ extending the mapping $x_i\to x_i$, $i=1,\ldots,m$.
Let $V$ be the $GL_m$-submodule of $K\langle X\rangle$ generated by the $h_j$.
Again, $\delta$ acts as a nilpotent linear operator on $V$. We fix a basis
$Z=\{z_1,\ldots,z_q\}$ of $V$. Consider the free $K[Y]$-module $M(Y,Z)$ with basis $Z$.
Proposition \ref{proposition 3} gives that $M(Y,Z)^{\delta}$
is a finitely generated $K[Y]^{\delta}$-module and the theorem of Weitzenb\"ock 
implies that $K[Y]^{\delta}$ is a finitely generated algebra. Since the algebra
$C_{nm}$ and the $C_{nm}$-module $T_{nm}$ are homomorphic images of the algebra
$K[Y]$ and the $K[Y]$-module $M(Y,Z)$, Corollary \ref{corollary of proposition 3}
gives that $K[Y]^{\delta}$ and $M(Y,Z)^{\delta}$ map on $C_{nm}^{\delta}$ and
$T_{nm}^{\delta}$, respectively. Hence $T_{nm}^{\delta}$ is a finitely generated
module of the finitely generated algebra $C_{nm}^{\delta}$ and, therefore, the algebra
$T_{nm}^{\delta}$ is finitely generated.
\end{proof}

\end{document}